\documentstyle[12pt,amscd,amssymb]{amsart}
\newtheorem{thm}{Theorem}[section]
\newtheorem{lemma}[thm]{Lemma}

\newtheorem{prop}[thm]{Proposition}
\newtheorem{rem}[thm]{Remark}

\newtheorem{defin}[thm]{Definition}

\begin{document}

\newcommand{\C}{{\mathbb C}}
\newcommand{\I}{{\o I}}
\newcommand{\CH}{{\mathbb H}_{\C}^2}
\renewcommand{\P}{{\mathbb P}}
\newcommand{\R}{{\mathbb R}}
\newcommand{\Z}{{\mathbb Z}}
\newcommand{\Os}{{\mathcal O}}
\newcommand{\zz}{{\mathcal Z}}
\newcommand{\z}{{\zz}^1(M)}
\newcommand{\M}{{\mathcal M}}
\newcommand{\Mpq}{{\mathcal M}_{(d_P,d_Q)}}
\newcommand{\Hpq}{{\mathcal H}_{(d_P,d_Q)}}
\newcommand{\la}{\langle}
\newcommand{\ra}{\rangle}
\newcommand{\h}{{\mathcal H}^3_{\R}}
\newcommand{\hH}{{\mathscr H}}
\newcommand{\hi}{h^{-1}}
\renewcommand{\o}{\operatorname}
\renewcommand{\Im}{\o{Im}}
\renewcommand{\Re}{\o{Re}}
\newcommand{\End}{\o{End}}
\newcommand{\Aut}{\o{Aut}}
\newcommand{\Hol}{\o{Hol}(E)}
\newcommand{\Her}{\o{Her}(E)}
\newcommand{\fc}{{\mathcal F}(E)} 
\newcommand{\fcu}{{\mathcal F}_u(E)} 
\newcommand{\fcl}{{\mathcal F}_l(E)} 
\newcommand{\ga}{{\mathcal G}(E)} 
\newcommand{\gau}{{\mathcal G}(E)} 
\newcommand{\gal}{{\mathcal G}_{\C}(E)} 
\newcommand{\Hig}{\o{Higgs}(V)}
\newcommand{\Ext}{\o{Ext}}
\newcommand{\Hom}{\o{Hom}}
\newcommand{\cHom}{{\mathcal H}\o{om}}
\newcommand{\Bpq}{{\mathcal B}_{(d_P,d_Q)}}
\newcommand{\Up}{\o{U}(p)}
\newcommand{\up}{{\frak u}(p)}
\newcommand{\PUp}{\o{PU}(p)}
\newcommand{\PUpp}{\o{PU}(p,p)}
\newcommand{\Uq}{\o{U}(q)}
\newcommand{\uq}{{\frak u}(q)}
\newcommand{\PUq}{\o{PU}(q)}
\newcommand{\Un}{\o{U}(n)}
\newcommand{\un}{{\frak u}(n)}
\newcommand{\PUn}{\o{PU}(n)}
\newcommand{\Uoo}{\o{U}(1,1)}
\newcommand{\uoo}{{\frak u}(1,1)}
\newcommand{\PUoo}{\o{PU}(1,1)}
\newcommand{\Ut}{\o{U}(2)}
\newcommand{\SUtt}{\o{SU}(2,2)}
\newcommand{\ut}{{\frak u}(2)}
\newcommand{\Uto}{\o{U}(2,1)}
\newcommand{\PUto}{\o{PU}(2,1)}
\newcommand{\uto}{{\frak u}(2,1)}
\newcommand{\Upq}{\o{U}(p,q)}
\newcommand{\Upp}{\o{U}(p,p)}
\newcommand{\Upo}{\o{U}(p,1)}
\newcommand{\upq}{{\frak u}(p,q)}
\newcommand{\upo}{{\frak u}(p,1)}
\newcommand{\PUpq}{\o{PU}(p,q)}
\newcommand{\PUpo}{\o{PU}(p,1)}
\newcommand{\Uo}{\o{U}(1)}
\newcommand{\uo}{{\frak u}(1)}
\newcommand{\GLn}{\o{GL}(n,\C)}
\newcommand{\gln}{{\frak gl}(n,\C)}
\newcommand{\PGLn}{\o{PGL}(n,\C)}
\newcommand{\PSLt}{\o{PSL}(3,\C)}
\newcommand{\Spf}{\o{Sp}(4,\R)}
\newcommand{\rk}{\o{rank}}
\newcommand{\Map}{\o{Map}}
\newcommand{\db}{\overline{\partial}}
\renewcommand{\d}{\partial}
\newcommand\ka{K\"ahler~}
\newcommand\tM{\tilde{M}}
\newcommand\ca{{\mathcal A}}
\newcommand\cb{{\mathcal B}}
\newcommand\cc{{\mathcal C}}
\renewcommand{\i}{{\o i}}
\renewcommand{\j}{{\o j}}
\renewcommand{\k}{{\o k}}
\newcommand\Iso{{\o{Iso}}}
\newcommand\Obj{{\o{Obj}}}
\newcommand\hk{{hyperk\"ahler~}}
\newcommand{\E}{{\cal E}}
\newcommand{\F}{{\cal F}}
\newcommand{\G}{{\cal G}}
\renewcommand{\H}{{\o H}}
\newcommand{\K}{{\cal K}}
\renewcommand{\L}{{\cal L}}
\newcommand{\N}{{\cal N}}
\newcommand{\U}{{\cal U}}
\newcommand{\V}{{\cal V}}
\newcommand{\W}{{\cal W}}
\renewcommand{\O}{{\cal O}}
\newcommand{\Q}{{\cal Q}}

\title
{The Moduli of Flat PU($p,p$)-Structures with Large Toledo Invariants}
\author{Eyal Markman
\and Eugene Z. Xia}
\address{ 
Department of Mathematics,
University of Massachusetts, Amherst, MA 01003-4515}
\email{markman@@math.umass.edu {\it (Markman)}, 
xia@@math.umass.edu {\it (Xia)}}
\date{\today} 
\thanks{E. Markman was partially supported by NSF grant number DMS-9802532}
\maketitle
\begin{abstract}
For a compact Riemann surface $X$ of genus $g > 1$, 
$\Hom(\pi_1(X), \PUpq)/\PUpq$ is the moduli space of 
flat $\PUpq$-connections on $X$.
There are two invariants, the Chern class $c$ and the Toledo invariant $\tau$
associated with each element in the moduli. 
The Toledo invariant is bounded in the range 
$-2min(p,q)(g-1) \le \tau \le 2min(p,q)(g-1)$.
This paper shows that 
the component,
associated with a fixed $\tau > 2(max(p,q)-1)(g-1)$ 
(resp. $\tau < -2(max(p,q)-1)(g-1)$)
and a fixed Chern class $c$, is connected (The restriction on $\tau$
implies $p=q$).
\end{abstract}


\section{Introduction and Results}
Let $X$ be a smooth projective curve over $\C$ with genus $g>1$.
Let $\PGLn$ and $\PUpq$ be $\GLn$ and $\Upq$ modulo their respective
centers.
The deformation space 
$$
\C {\cal N}_B = \Hom^+(\pi_1(X), \PGLn)/\PGLn
$$
is the space of equivalence classes of 
semi-simple $\PGLn$-representations of the fundamental
group $\pi_1(X)$.
This is the $\PGLn$-Betti moduli space on $X$. 

Since $\PUpq \subset \PGLn$, $\C {\cal N}_B$ contains the space
$$
{\cal N}_B(p,q) = \Hom^+(\pi_1(X), \PUpq)/\PUpq.
$$
The space ${\cal N}_B(p,q)$ will 
be referred to as the $\PUpq$-Betti moduli space.
One may assume $p \ge q > 0$ without loss of generality.

The Betti moduli spaces are of great interest in
geometric topology and uniformization.
Goldman analyzed ${\cal N}_B(1,1)$
and determined the number of its connected
components \cite{Go1}.
A theorem of Corlette, Donaldson, Hitchin and Simpson 
gives a homeomorphism of $\C {\cal N}_B$ to
two other moduli spaces---the 
$\PGLn$-de Rham and the $\PGLn$-Dolbeault moduli spaces, respectively
\cite{Co1,Do0, Hi1, Si1}.  The Dolbeault moduli spaces are moduli
of semi-stable Higgs bundles.  Hitchin subsequently considered 
${\cal N}_B(1,1)$ from the Higgs bundle perspective
and determined its topology \cite{Hi1}.
The cases when the structure groups being $\Upo$ and $\PUto$
are treated in \cite{Xi3, Xi4}.  Gothen obtained partial 
results for the structure groups $\SUtt$ and $\Spf$ \cite{Go2}.
Other related results have been obtained in \cite{Xi1, Xi2}.

Each element in ${\cal N}_B(p,q)$ is associated with a Chern class
$c$ and a Toledo invariant
$\tau$ which is bounded as \cite{Do1,To1,To2}
$$
-2 q (g-1) \le \tau \le 2 q (g-1).
$$
The main result presented here is the following:
\begin{thm} \label{thm:main}
The locus in ${\cal N}_B(p,q)$, associated with  
a fixed $\tau > 2(p-1)(g-1)$ (resp. $\tau < -2(p-1)(g-1)$)
and a fixed Chern class $c$, is connected.
\end{thm}
\begin{rem} \label{rem:p=q}
The hypothesis $\tau>2(p-1)(g-1)$ and
the fact that $0 \le \tau \le 2 q (g-1)$ imply $p=q$.
\end{rem}

Section 
\ref{sec-background} reviews the homeomorphism between the
$\GLn$-Betti space and the moduli space of Higgs bundles. 
Section \ref{sec-Upp-higgs-bundles} recalls the 
characterization of $\Upq$-Higgs bundles and their Toledo invariant. 
Section \ref{sec-Upp-Hodge-bundles} concerns
the $\C^*$-action on the moduli space of
$\Upp$-Higgs bundles and introduces
the locus of Binary hodge-bundles, a distinguished component
of the $\C^*$-invariant locus. The main result of Section
\ref{sec-Upp-Hodge-bundles} (Proposition \ref{prop:binary}) implies that any point 
in the components (with $\tau > 2(p-1)(g-1)$) of the moduli can 
be deformed to a binary hodge bundle.  
In Section \ref{sec:binary-is-connected}, we prove that the locus of binary hodge bundles 
is irreducible (Proposition~\ref{prop:irreducible}).  
Theorem~\ref{thm:main} follows from Propositions~\ref{prop:binary}
and \ref{prop:irreducible}.

After the completion of this paper, S. Bradlow, O. Garcia-Prada
and P. Gothen announced that all moduli spaces of flat $\PUpq$-structures
with fixed Chern class and Toledo invariant, are connected \cite{Br}.

\

\centerline{\sc Acknowledgments}
We thank the referee for suggestions for improvement.  We also thank
Steven Bradlow, Oscar Garcia-Prada and Peter Gothen for pointing out a 
gap in a previous version.  Eugene Z. Xia thanks the National
Center for Theoretical Sciences, Taiwan for hospitality.

\section{The $\GLn$-Higgs Bundles}
\label{sec-background}

Let $\Gamma$ be the central extension
$$
1 \longrightarrow \Z \longrightarrow \Gamma \longrightarrow \pi_1(X) 
\longrightarrow 1
$$
as in \cite{At1, Hi1}.
Each $\rho \in \Hom(\Gamma, \GLn)$ acts on $\C^n$
via the standard representation of $\GLn$.  
The representation $\rho$ is called
reducible (resp. irreducible),
if its action on $\C^n$ is reducible (resp. irreducible).  A representation
$\rho$ is called semi-simple, 
if it is a direct sum of irreducible representations.
\begin{defin}
$$
\C {\cal M}_B = \{\sigma \in \Hom(\Gamma, \GLn) : \sigma
\mbox{ is semi-simple}\}/\GLn.
$$
$$
{\cal M}_B(p,q) = \{\sigma \in \Hom(\Gamma, \Upq) : \sigma
\mbox{ is semi-simple}\}/\Upq.
$$
\end{defin}
It is immediate that 
$$
\C {\cal N}_B  = \C {\cal M}_B/\Hom(\pi_1(X), \C^*)
$$
$$
{\cal N}_B(p,q) = {\cal M}_B(p,q)/\Hom(\pi_1(X), \Uo).
$$
Therefore counting the components of ${\cal N}_B(p,q)$ is the same
as counting the components of ${\cal M}_B(p,q)$.

Let $E$ be a rank-$(p+q)$ complex
vector bundle over $X$ with
$0 \le \deg(E) < p + q$.  Denote by $\Omega$ the 
canonical bundle on $X$.  A holomorphic structure
$\db$ on $E$ induces holomorphic structures on the bundles
$\End(E)$ and $\End(E) \otimes \Omega$.
A Higgs bundle is a pair $(E_{\db}, \Phi)$,
where $\db$ is a holomorphic structure on $E$ and 
$\Phi \in \H^0(X, \End(E_{\db}) \otimes \Omega)$.  Such 
a $\Phi$ is called a Higgs field.
We denote the holomorphic bundle $E_{\db}$ by $V$.

Define the slope of a vector bundle $V$ to be
$$
s(V) = \deg(V)/\rk(V).
$$
For a fixed $\Phi$, a holomorphic sub-bundle
$W \subset V$ is said to be $\Phi$-invariant if
$\Phi(W) \subset W \otimes \Omega$.
A Higgs bundle $(V,\Phi)$ is stable (semi-stable) if $W \subset V$ being
$\Phi$-invariant implies
$$
s(W) < (\le) s(V).
$$
A Higgs bundle is called poly-stable if it is a direct
sum of stable Higgs bundles of the same slope \cite{Hi1,Si2}.

The Dolbeault moduli space
$\C {\cal M}$ is the coarse moduli space 
of semi-stable rank-$(p+q)$ Higgs bundles
on $X$ \cite{Hi1, Hi2, Ni1, Si2}.  The closed points of $\C {\cal M}$
parameterize the $S$-equivalent classes of semi-stable Higgs bundles.
Moreover every $S$-equivalence class has a poly-stable representative.

\section{The $\Upp$-Higgs Bundles}
\label{sec-Upp-higgs-bundles}
\begin{defin} \label{def:M}
Let $\M$ be the subset of  $\C {\cal M}$ consisting
of equivalent classes of Higgs bundles, whose poly-stable 
representative $(V, \Phi)$ satisfies the following two conditions:
\begin{enumerate}
\item $V$ is a direct sum:
$$
V = V_P \oplus V_Q,
$$
where $V_P, V_Q$ are of ranks $p$ and $q$, respectively.
\item The Higgs field decomposes into two maps:
$$
\Phi_1 : V_P \longrightarrow V_Q \otimes \Omega,
$$
$$
\Phi_2 : V_Q \longrightarrow V_P \otimes \Omega.
$$
\end{enumerate}
\end{defin}
Hence each
$(V,\Phi) = (V_P \oplus V_Q, \Phi) \in \M$
is associated with two invariants, $d_P = \deg(V_P)$ and $d_Q = \deg(V_Q)$.
The Toledo invariant is defined to be 
$$
\tau = 2 \frac{\deg(V_P \otimes V_Q^*)}{p+q} 
= 2 \frac{q d_P - p d_Q}{p+q},
$$
and the Chern class is $d_P + d_Q$.
The subset of ${\cal M}$, consisting of classes with fixed $d_P$ and $d_Q$, 
is denoted by $\Mpq$.

\begin{rem}
\label{rem-ambiguous-notation}
{\rm
When $p=q$, the labeling of the summands $V_P$ and $V_Q$ may seem ambiguous. 
However, we will introduce assumption
(\ref{large-tau}) below, which implies that $d_P>d_Q$. 
Thus, when $p=q$, $V_P$ is
distinguished as the summand of larger degree. 
}
\end{rem}

\begin{thm} \label{thm:summary}
\mbox{}

\begin{enumerate}
\item The moduli spaces $\C {\cal M}_{B}$ and
$\C {\cal M}$ are homeomorphic.
\item The reducible representations in $\C {\cal M}_{B}$ correspond to
the poly(semi)-stable, but not stable, points.
\item 
\label{thm-item-fixed-locus-of-an-involution}
The subspace ${\cal M}_{B}$ is homeomorphic to $\M$.
\item 
\label{thm-item-symmetry-of-invariants}
The space $\Mpq$ is homeomorphic to
${\cal M}_{(-d_P,-d_Q)}$.
\end{enumerate}
\end{thm}
\begin{pf}
The proof of (1) and (2) can be found in \cite{Co1} (the main idea 
is present in
\cite{Do0, Hi1} and the most general version of this celebrated result
is in \cite{Si1}).  For (2) and (3), see \cite{Si1, Xi3}.  See also the
last few sections of \cite{Si2} for general real forms.
\end{pf}

Part (\ref{thm-item-fixed-locus-of-an-involution}) characterizes isomorphism
classes of semi-stable
$\Upq$-Higgs bundles as those, which are fixed under the involution
$(V,\Phi)\mapsto (V,-\Phi)$ and the involution of $V$ conjugating
$\Phi$ to $-\Phi$ has eigenvalue $-1$ with multiplicity $q$. 
In particular, ${\cal M}$ is a closed subset of the 
quasi-projective variety $\C {\cal M}$, and we may endow ${\cal M}$ 
with the reduced induced subscheme structure. 

By Theorem~\ref{thm:summary} (\ref{thm-item-symmetry-of-invariants}), 
we assume, for the rest of the paper that 
$$0 \le \tau \le 2 q (g-1).$$

\begin{lemma} \label{prop-Phi-1-is-generically-surjective}
Suppose $\tau > 2(p-1) (g-1)$.  Then 
$(V_P \oplus V_Q,(\Phi_1, \Phi_2)) \in \Mpq$
implies
$$
\Phi_1 : V_P \longrightarrow V_Q \otimes \Omega
$$
is generically surjective.
\end{lemma}
Lemma~\ref{prop-Phi-1-is-generically-surjective}
was first obtained for integer $\tau$ by Gothen 
\cite{Go2}.
\begin{pf}
Since $(V_P \oplus V_Q,\Phi)$ is semi-stable and $\tau$ is positive,
$\Phi_1 \not\equiv 0$.
Suppose $\Phi_1$ is not generically surjective with a non-trivial kernel $V_1$.  
We want to produce a destabilizing Higgs subbundle, namely, 
$(V_P \oplus W_1 \otimes \Omega^{-1}, \Phi)$, where $W_1$ is as in
the following canonical factorization of $\Phi_1$:
$$
\begin{array}{ccccccccc}
0 & \longrightarrow & V_1 & \stackrel{f_1}{\longrightarrow} &
V_P & \stackrel{f_2}{\longrightarrow} &
V_2 & \longrightarrow & 0\\
  & & & & \Phi_1 \Big\downarrow & &
\varphi \Big\downarrow & & \\
0 & \longleftarrow & W_2 & \stackrel{g_2}{\longleftarrow}
& V_Q \otimes \Omega & \stackrel{g_1}{\longleftarrow} &
W_1 & \longleftarrow & 0
\end{array}.
$$
In the above diagram, the rows are exact, $\rk(V_2) = \rk(W_1)$ and
$\varphi$ has full rank at a generic point of $X$.
Let $d_i = \deg(V_i)$, $r_i = \rk(V_i)$.  
Since $V_1$ is $\Phi$-invariant, semi-stability 
implies 
$$
\frac{d_1}{r_1} \le \frac{d_P+d_Q}{2p}.  
$$
Since $\varphi$ has full rank generically,
$\deg(W_1) \ge d_2 = d_P - d_1$.
Hence, 
\begin{eqnarray*}
\deg(V_P \oplus W_1 \otimes \Omega^{-1}) & = & d_P + \deg(W_1) - 2 r_2 (g-1)\\
& \ge & 2d_P - \frac{r_1(d_P+d_Q)}{2p} - 2 r_2 (g-1).
\end{eqnarray*}
This implies
\begin{eqnarray*}
s(V_P \oplus W_1 \otimes \Omega^{-1}) - s(V) & \ge & 
\frac{2d_P - \frac{r_1(d_P+d_Q)}{2p} - 2 r_2 (g-1)}{p+r_2} - 
\frac{d_P + d_Q}{2p}.
\end{eqnarray*}
This, together with the facts 
$r_1 + r_2 = p$
and $\tau = d_P - d_Q$, gives us
$$
s(V_P \oplus W_1 \otimes \Omega^{-1}) - s(V) \ge 
\frac{\tau - 2 r_2 (g-1)}{p+r_2}.
$$
Since $r_1 \ge 1$, we have $r_2 < p$.  The assumption $\tau > 2 (p-1) (g-1)$
then implies that $V_P \oplus W_1 \otimes \Omega^{-1}$
destabilizes $V$.
\end{pf}
For the rest of the paper, we assume 
\begin{equation} \label{large-tau}
\tau > 2 (p-1) (g-1).
\end{equation}

\section{The $\Upp$-Hodge Bundles}
\label{sec-Upp-Hodge-bundles}

\begin{defin}[See \cite{Si2}]
A Hodge bundle on $X$
is a Higgs bundle $(E,\Phi)$ of the following form:
$$
E = \bigoplus_{i=0}^k E^i
$$
and $\Phi = (\phi_k,...,\phi_1)$ with:
$$
\phi_i : E^i \longrightarrow E^{i-1} \otimes \Omega.
$$
The integer $k$ is called the {\em length} of the Hodge bundle
$(E,\Phi)$.
\end{defin}

There is a $\C^*$-action on $\C{\cal M}$:
\begin{eqnarray*}
\C^* \times \C {\cal M} & \longrightarrow & \C {\cal M}
\\
(t,E,\Phi) & \mapsto & (E,t\Phi).
\end{eqnarray*}
The space of equivalence classes of semi-stable Hodge
bundles is the set of fixed points of this action 
(See Lemma 4.1 of \cite{Si2}).

\begin{rem}
{\rm
If $(E,\Phi)$ is a stable Higgs bundle, which is a fixed point of the 
$\C^*$-action, then it admits a {\em unique} decomposition of $E$,  
realizing it as a Hodge bundle (Lemma 4.1 in \cite{Si2}). 
Hence, the decomposition is canonical if $(E,\Phi)$ is poly-stable.
}
\end{rem}

A stable Hodge bundle $(E,\Phi)$ 
admits a unique realization as a $\Upq$-Higgs bundle, 
for a unique pair of non-negative integers $(p,q)$. 
Let $E=\oplus_{i=0}^k E^i$ be the unique 
decomposition of $(E,\Phi)$ and
$\phi_i : E^i \rightarrow E^{i-1}$ the corresponding
decomposition of the Higgs field.
Let $E^{{\rm odd}}$ be the direct sum of the summands with odd
index. We define $E^{{\rm even}}$ similarly.
The stability of $(E,\Phi)$ implies, that there exists
at most one automorphism $f$ of $E$, up to a $\C^*$-factor, 
which conjugates $\Phi$ to $-\Phi$. Such an
automorphism $f$ is given by multiplying $E^{{\rm odd}}$ by 
$a\in \C^*$ and multiplying $E^{{\rm even}}$ by $-a$. 
Set $V_P=E^{{\rm odd}}$ and $V_Q=E^{{\rm even}}$.
Then $(V_P\oplus V_Q,\Phi)$ is a $\Upq$-Higgs bundle,
with $p={\rm rank}(V_P)$ and $q={\rm rank}(V_Q)$. 
Recall our convention, that if $p=q$, 
then $V_P$ is the summand of larger degree
(Remark \ref{rem-ambiguous-notation}). 
This follows from assumption (\ref{large-tau})  
by Lemma \ref{prop-Phi-1-is-generically-surjective}.
Otherwise,  
$\deg(E^{{\rm even}})>\deg(E^{{\rm odd}})$
and Lemma \ref{prop-Phi-1-is-generically-surjective} implies that 
$\Phi : E^{{\rm even}}\rightarrow E^{{\rm odd}}\otimes \Omega$ is injective.
This would contradict the fact that $E^0$ is in the kernel of $\Phi$.

In Section \ref{sec-C-star-action} 
we describe the $\C^*$-action on the tangent space of a Hodge bundle.
In Section \ref{sec-binary-hodge-bundles} 
we introduce the locus of Binary Hodge bundles
and characterize it in terms of the $\C^*$-action. 

\subsection{The infinitesimal $\C^*$-action} 
\label{sec-C-star-action}

\medskip
The infinitesimal deformations of a Higgs bundle $(E,\Phi)$ are calculated
by the first cohomology of the complex $K_\bullet$ below
\[
\End(E) \ \ \ \stackrel{ad_\Phi}{\longrightarrow} \ \ \ \End(E)\otimes\Omega
\]
(in degrees $0$ and $1$). When $(E,\Phi)$ is a Hodge bundle, $\H^1(K_\bullet)$
decomposes into weight spaces of the natural $\C^*$-action.
Next, we analyze this decomposition.
Let
\begin{eqnarray*}
\mu \ : \  \C^* & \rightarrow & \Aut(E^i) \ \ \ \mbox{and}
\\
\alpha \ : \  \C^* & \rightarrow & \Aut(\Omega)
\end{eqnarray*}
be the representations with weights $i$ and $1$ respectively.
We denote by $\mu$ also its natural extension to tensor products of the $E^i$
and by $\alpha$ the action on $\Hom(E^i,E^j)\otimes\Omega$ via
$1\otimes\alpha$. Set
\[
\rho \ := \ \mu\cdot \alpha.
\]
\noindent
Then $\mu$ and $\rho$ both have weight $j-i$ on $\Hom(E^i,E^j)$.
On $\Hom(E^i,E^j\otimes\Omega)$, $\mu$ has weight $j-i$ and
$\rho$ has weight $j-i+1$.
Observe that the Higgs field $\Phi$, of a Hodge bundle,
has weight $-1$ with respect to $\mu$ and
it is $\rho$-invariant. Consequently, the differential of the complex
$K_\bullet$ is $\rho$-invariant and the complex decomposes as a direct sum
\begin{equation} \label{complex-K-dot}
K_\bullet \ \ = \ \ \bigoplus_{w=-k}^{k+1} K_\bullet^w,
\end{equation}
where $K_\bullet^w$ is the complex
\[
\bigoplus_{i=\max\{0,-w\}}^{\min\{k,k-w\}}\Hom(E^i,E^{i+w}) \ \ \
\stackrel{ad_\Phi}{\longrightarrow}  \ \ \
\bigoplus_{i=\max\{0,1-w\}}^{\min\{k,k-w+1\}}\Hom(E^i,E^{i+w-1}\otimes\Omega).
\]
For example, $K_\bullet^{-k}$ is the complex supported in degree zero
by the vector bundle
$\Hom(E^k,E^0)$.
The complex $K_\bullet^{1-k}$ is
\begin{eqnarray}
\label{eq-complex-K-1-k}
\Hom(E^{k-1},E^0)\bigoplus \Hom(E^k,E^1) & \stackrel{ad_\Phi}{\longrightarrow}  &
\Hom(E^k,E^0\otimes\Omega),
\\
\nonumber
(a,b) \hspace{12ex} & \mapsto & \phi_1\circ b - (a\otimes 1)\circ \phi_k.
\end{eqnarray}

\begin{lemma}
The representation $\rho:\C^* \rightarrow \Aut(\H^1(K_\bullet))$
is the infinitesimal
$\C^*$-action, on the tangent space at the fixed point $(E,\Phi)$,
arising from the $\C^*$-action $(F,\varphi)\mapsto  (F,t\varphi)$
on the moduli of Higgs bundles.
\end{lemma}
\begin{pf}
Denote by $K_{\bullet,t}$ the complex corresponding to deformations
of the Higgs pair $(E,t^{-1}\Phi)$. There are {\em two}
natural isomorphisms of complexes
\begin{eqnarray}
\label{eq-relating-two-cohomological-identifications}
\mu_t \ : \ \H^1(K_\bullet) & \stackrel{\cong}{\longrightarrow} 
& \H^1(K_{\bullet,t})
\ \ \ \mbox{and}
\\
\nonumber
(1,t^{-1}) \ : \ \H^1(K_\bullet) & \stackrel{\cong}{\longrightarrow}
 & \H^1(K_{\bullet,t})
\end{eqnarray}
The first is the evaluation of the representation $\mu$ at $t$.
The second is the identity on the vector bundle in degree $0$ and
multiplication by $t^{-1}$ on the vector bundle in degree $1$.
The automorphism $\mu_t$ of $E$ sends $\Phi$ to $t^{-1}\Phi$.
Hence, $(E,\Phi)$ and
$(E,t^{-1}\Phi)$ are equivalent Hodge bundles.
In other words, Hodge bundles are $\C^*$-invariant.
The isomorphism
(\ref{eq-relating-two-cohomological-identifications})
is the natural identification of the two cohomological calculations of the same
tangent space.
The automorphism $\rho_t$ of $K_\bullet$ is the composition
\[
\rho_t \ \ = \ \ (1,t^{-1})^{-1}\circ \mu_t.
\]
The lemma follows.
\end{pf}

Denote by $K^{\rm even}_\bullet$
the sum of the even weight sub-complexes of $K_\bullet$.
\begin{lemma} \label{lemma-infinitesimal-u-p-q-deformations}
Assume that $(E,\Phi)$ is stable. Then
the infinitesimal $U(p,q)$-deformations of a Hodge bundle $(E,\Phi)$
are parameterized by $\H^1(K^{\rm even}_\bullet)$.
\end{lemma}
\begin{pf}
This is the infinitesimal counterpart of the characterization of
the $U(p,q)$ moduli space as the fixed point set of the involution
\[
(E,\Phi) \ \ \mapsto \ \ (E,-\Phi).
\]
\end{pf}

\subsection{Binary Hodge bundles}
\label{sec-binary-hodge-bundles}

The  $\C^*$-action on $\C {\cal M}$ 
preserves $\Mpq$ \cite{Si2,Xi3}.  Let 
$\C {\cal H} \subset \C {\cal M}$ and $\Hpq \subset \Mpq$ be
the corresponding sets of Hodge bundles.  
\begin{defin} \label{defin:hodge}
Suppose $(E,\Phi) = (V_P \oplus V_Q, (\Phi_1, \Phi_2)) \in  \Mpq$.
Then $(V_P \oplus V_Q, (\Phi_1, \Phi_2))$ is a binary Hodge bundle
if $\Phi_2 \equiv 0$.  That is a binary hodge bundle is a $\Upq$-Hodge
bundle of length $1$.
Denote by $\Bpq$ the set of binary Hodge bundles.
\end{defin}

$\Bpq$ is a closed subset of $\Mpq$. 
It may be endowed with the reduced induced subscheme structure. 
The closedness of $\Bpq$ follows from 
Lemma \ref{lemma:deform} below. Lemma \ref{lemma:deform} characterizes 
$\Bpq$ as the union of all connected components of the fixed locus 
of $\Mpq$, which are the minimal with respect to a 
partial ordering by the weight invariant (\ref{eq-weight-invariant}). 
We show later that $\Bpq$ is irreducible (Proposition \ref{prop:irreducible}).

The rest of this section is dedicated to proving the following proposition.

\begin{prop} \label{prop:binary}
Under the assumption (\ref{large-tau}), $\Bpq$
intersects every connected component of $\Mpq$.
\end{prop}

The proposition follows from Lemmas \ref{lem:deform} 
and \ref{lemma:deform}. 
We will follow the general outline of Simpson's proof of the connectedness of 
the moduli space $\C {\cal M}$ (see \cite{Si4}). Simpson introduced
an algebraic version of a standard Morse theoretic technique.  The space
$\C {\cal M}$, as well as its subvariety $\Mpq$, are 
quasi-projective varieties. They fit into the following set-up.
Let $Y$ be a quasi-projective variety upon which $\C^*$ acts 
algebraically. Assume that $L$ is a very ample line-bundle on $Y$ 
with a linearization of the action. Then there is a 
finite dimensional invariant subspace $V\subset H^0(Y,L)$, giving rise to  a 
$\C^*$-equivariant embedding $Y\hookrightarrow \P V^*$.
Denote by $V=\oplus V_\alpha$ the decomposition into weight subspaces,
so that $t\in \C^*$ acts by $t^\alpha$ on $V_\alpha$. 
Let $Z$ be the closure of $Y$ in $\P V^*$. 
Denote by $Z_\beta$ the intersection $Z\cap \P V_\beta^*$. 
The fixed locus of $Z$ is the union of the loci $Z_\beta$. In particular, 
a connected component $Z'$, of the fixed locus of $Z$, 
comes with an invariant; a weight 
\begin{equation}
\label{eq-weight-invariant}
\beta(Z')
\end{equation}
such that $Z'$ is contained in $\P V_{\beta(Z')}^*$.
Assume $y\in Y$ is not a fixed point. Then the closure of the 
$\C^*$-orbit of $y$ in $Z$ has two fixed points
$y_0:=\lim_{t\rightarrow 0}ty$ and 
$y_\infty := \lim_{t\rightarrow \infty}ty$.
Moreover, the invariants of $y_0$ and $y_\infty$ satisfy the strict inequality
$\beta(y_0)<\beta(y_\infty)$. Thus, connected components of the
fixed locus are partially ordered. Moreover, if $y$ is not a fixed point, 
taking the limit $\lim_{t\rightarrow 0}ty$ amounts to
{\em flowing down to a lower connected component of the fixed locus of $Z$}. 

Assume, furthermore, 
that $\lim_{t \rightarrow 0} tx$ exists in $Y$ for all $x \in Y$. 
Then the process of flowing down can be used to study the connectedness of
$Y$. More precisely, we have the following lemma, 
which is a trivial variation on Lemma 11.8 in \cite{Si4}:

\begin{lemma} [Lemma 11.8 in \cite{Si4}] \label{lem:deform}
Suppose $Y$ is a quasi-projective variety, upon which $\C^*$ acts 
algebraically and $\lim_{t \rightarrow 0} tx$ exists in $Y$ for all $x \in Y$ as above. 
Suppose $U \subset Y$ is a subset of the fixed point set of $\C^*$, and
suppose that for any fixed point $x \not\in U$, there exists $y \neq x$
such that $\lim_{t \rightarrow \infty} ty = x$.  Then $U$ intersects every
connected component of $Y$. 
\end{lemma}

\begin{pf}
Suppose $Y'$ is a connected component of $Y$ not intersecting $U$.
Let $\beta$ be the smallest integer such that $Y'_\beta$ is non-empty.
Choose $x\in Y'_\beta$. 
By hypothesis, there exists $y\neq x$ in $Y'$, 
such that $\lim_{t\rightarrow \infty}ty=x$. On the other hand, 
$z':=\lim_{t\rightarrow 0}ty$ is also in $Y'$, say in $Y'_\alpha$. 
Since $y$ is not a fixed point, $\alpha<\beta$. This contradicts the 
minimality of $\beta$. 
\end{pf}

It is known that $\lim_{t \rightarrow 0} (E, t\Phi)$
always exists in $\C {\cal M}$ and is a Hodge bundle \cite{Si2}. 
The same holds for $\Mpq$ since $\Mpq$ is
closed in $\C \M$ and is $\C^*$-invariant.
Proposition \ref{prop:binary} follows from 
Lemma \ref{lem:deform}, with $Y=\Mpq$ and $U=\Bpq$,
and Lemma \ref{lemma:deform}. The latter is 
the analogue of Lemma 11.9 of \cite{Si4}.

%


\begin{lemma} \label{lemma:deform}
Suppose $(E,\Phi) \in \Hpq$ is  poly-stable with length $k>1$.  
Then there exists $(F,\Psi) \in \Mpq$ satisfying 
\begin{eqnarray*}
\lim_{t \rightarrow \infty} (F,t\Psi) & \cong & (E,\Phi) \ \ \ \mbox{and}
\\
(F,\Psi) & \not\cong & (E,\Phi).
\end{eqnarray*}
\end{lemma}

\begin{pf}
Assume first that $(E,\Phi)$ is stable. Then it is a smooth point of $\Mpq$.
Lemma \ref{lemma-existence-of-negative-weights} below implies that 
negative weights occur in the weight decomposition of the tangent space
to $\Mpq$ at $(E,\Phi)$. 
Take a $\C^*$-orbit $R$ in $\Mpq$, with $(E,\Phi)$ in its closure, 
such that the tangent line to $R$ at $(E,\Phi)$ has negative weight. 
Then any point $(F,\Psi)$ in $R$ would satisfy the conditions of the Lemma.

Next, we reduce the proof to the stable case. 
Suppose $(E,\Phi)=(V_P\oplus V_Q,\Phi)$
is a poly-stable representative, of length $k>1$, 
of an equivalence class in $\Hpq$. 
Then $(E,\Phi)$ is of the form
$(V'_P\oplus V'_Q,\Phi') \oplus (V''_P\oplus V''_Q,\Phi'')$, 
where $(V'_P\oplus V'_Q,\Phi')$ is a stable Hodge
bundle of length $k>1$. 
Lemma \ref{prop-Phi-1-is-generically-surjective} implies that 
$\Phi_1$ is generically an isomorphism. Hence, the same holds for 
$\Phi'_1$ and $\Phi''_1$. Lemma \ref{lemma-existence-of-negative-weights} 
below implies the 
existence of a pair $(F',\Psi')$, not equivalent to 
$(V'_P\oplus V'_Q,\Phi')$, such that
$\lim_{t\rightarrow \infty}(F',t\Psi')=(V'_P\oplus V'_Q,\Phi')$. 
Take $(F,\Psi):=(F',\Psi')\oplus (V''_P\oplus V''_Q,\Phi'')$.
Then 
$$
\lim_{t\rightarrow \infty} (F,t\Psi) = 
\lim_{t\rightarrow \infty} (F',t \Psi')\oplus (V''_P\oplus V''_Q,\Phi'') = (E,\Phi). 
$$
Moreover,
the graded objects $gr((E,\Phi))$ and $gr((F,\Psi))$
are not isomorphic. Hence, $(E,\Phi)$ and $(F,\Psi)$
are not equivalent.
\end{pf}

Let $(E,\Phi)$ be a stable hodge bundle corresponding to a 
$\Upq$-bundle $(V_P\oplus V_Q,\Phi)$. Since 
$(E,\Phi)$ is fixed by the $\C^*$-action, 
the tangent space $T_{(E,\Phi)}\Mpq$ is a representation of $\C^*$.
Denote by $[T_{(E,\Phi)}\Mpq]^-$ the sum of negative weight spaces
of this representation. Lemma \ref{lemma-existence-of-negative-weights}
below is used in the proof of Lemma \ref{lemma:deform}. 

\begin{lemma} \label{lemma-existence-of-negative-weights}
Let $(E,\Phi)$ be a stable hodge bundle corresponding to a 
$\Upp$-bundle $(V_P\oplus V_Q,\Phi)$. Assume that
$\Phi_1:V_P\rightarrow V_Q\otimes\Omega$ is an injective homomorphism. 
Assume further that the length $k$ of its Hodge decomposition is $\geq 2$. 
Then
$
[T_{(E,\Phi)}\Mpq]^-
$ 
does not vanish.
\end{lemma}

\begin{pf} Note that we exclude the
binary Hodge bundles which correspond to the case of $k=1$.
There are two cases, namely when $k$ is even and when $k$ is odd.

\noindent {\em Case 1}: $k > 2$ is odd. 

\noindent 
It suffices to prove that $\H^1(K_\bullet)^{1-k}$  does not vanish.
Our assumptions imply that
$\phi_{2i+1}:E^{2i+1}\rightarrow E^{2i}$ is generically an isomorphism.
In particular, $r_{2i+1}=r_{2i}$, where $r_i$ is the rank of $E^i$.
Stability of the Hodge bundle $(E,\Phi)$ implies:

\[
s(E^0\oplus E^1) \ \ \ < \ \ \ s(E) \ \ \ < \ \ \ s(E^{k-1}\oplus E^k).
\]

\begin{lemma}
\label{lemma-one-of-two-inequalities-holds}
At least one of the two cases holds:
\begin{eqnarray}
\label{eq-even-case}
s(E^0) & < & s(E^{k-1})  \ \ \ \mbox{or}
\\
\label{eq-odd-case}
s(E^1) & < & s(E^k).
\end{eqnarray}
\end{lemma}
\begin{pf}
Assume otherwise. Then $s(E^0) \geq s(E^{k-1})$ and
$s(E^1) \geq s(E^k)$. Using the fact that $r_0=r_1$ and $r_{k-1}=r_k$,
we get:
\[
s(E^0\oplus E^1) = \frac{s(E^0)+s(E^1)}{2} \geq
\frac{s(E^{k-1})+s(E^k)}{2} = s(E^{k-1}\oplus E^k).
\]
This contradicts the stability of $(E,\Phi)$.
\end{pf} 
We will use Lemma \ref{lemma-one-of-two-inequalities-holds} 
to prove that $\H^1(K_\bullet)^{1-k}$  does not vanish.
The $1-k$ weight space is equal to the first cohomology
$\H^1(K^{1-k}_\bullet)$ of the complex (\ref{eq-complex-K-1-k}).
$\H^i(K^{1-k}_\bullet)$ vanishes for $i<0$ and $i>2$. We get the
inequality
\[
\dim \H^1(K_\bullet)^{1-k} \ \ \ \geq \ \ \ -\chi(K^{1-k}_\bullet).
\]
Consequently, it suffices to prove that the Euler characteristic
$\chi(K^{1-k}_\bullet)$ is negative

\begin{equation}
\label{eq-chi-K-1-k-is-negative}
\chi(\Hom(E^{k-1},E^0))+\chi(\Hom(E^k,E^1))
-\chi(\Hom(E^k,E^0\otimes \Omega)) \ \ \ < \ \ 0.
\end{equation}

Consider first the case $s(E^1)<s(E^k)$. The Euler characteristic
$\chi(\Hom(E^k,E^1))$ is $r_1r_k[s(E^1)-s(E^k)+1-g]$.
It follows that $\chi(\Hom(E^k,E^1))$ is negative.
Composition with $\phi_k$
\[
\Hom(E^{k-1},E^0) \ \ \ \stackrel{\phi_k}{\longrightarrow} \ \ \
\Hom(E^k,E^0\otimes \Omega)
\]
is generically an isomorphism by our assumption on $\Phi_1$.
Consequently,
the difference of their Euler characteristics is negative
\begin{eqnarray*}
\chi(\Hom(E^{k-1},E^0))-\chi(\Hom(E^k,E^0\otimes \Omega)) & = &
\\
r_0r_k[s(\Hom(E^{k-1},E^0))-s(\Hom(E^k,E^0\otimes \Omega))] & < & 0.
\end{eqnarray*}
Equation (\ref{eq-chi-K-1-k-is-negative}) follows.

The case $s(E^0)<s(E^{k-1})$ is similar. Use composition with
$\phi_1$ instead.

\noindent {\em Case 2}: $k$ is even.

\noindent Then the complex $K_\bullet^{-k}$ is simply the vector
bundle $\Hom(E^k,E^0)$ supported in degree $0$.  Stability implies
\[
s(E^0) \ \ \ < \ \ \ s(E) \ \ \ < \ \ \ s(E^k).
\]
If follows that
$$
\dim \H^1(K_\bullet)^{-k} \ge -\chi(K_\bullet^{-k}) = 
-\chi(\Hom(E^k,E^0)) > 0.
$$
This completes the proof of Lemma 
\ref{lemma-existence-of-negative-weights}.
\end{pf}


\section{Connectedness of the locus of $\Upp$-Binary Hodge Bundles} 
\label{sec:binary-is-connected}
In this section, we show that $\Bpq$ is irreducible for fixed $d_P$ and
$d_Q$, thus, proving Theorem \ref{thm:main}.

\begin{prop} \label{prop:irreducible}
If $\tau>2(p-1)(g-1)$, then $\Bpq$ is irreducible.
\end{prop}
Recall that the assumption of Lemma~\ref{prop-Phi-1-is-generically-surjective} implies $p=q$
(See Remark \ref{rem:p=q}).
Key to the proof of Proposition~\ref{prop:irreducible} is the alternative 
description of $\Bpq$ provided by Lemma \ref{prop-Phi-1-is-generically-surjective}.
Suppose $(V_P \oplus V_Q, \Phi)$ is a binary Hodge bundle.
Then we obtain the length $p(2g-2)+d_Q-d_P$ quotient sheaf 
$(V_Q, f : V_Q \rightarrow V_Q/(\Phi_1(V_P)\otimes \Omega^{-1}))$.
Conversely, for each
length $p(2g-2)+d_Q-d_P$ quotient sheaf $(E,f : E  \rightarrow F)$,
we obtain a binary Hodge bundle $(\ker(f)\otimes\Omega \oplus E,\Phi_1)$,
where 
$$
\Phi_1 : \ker(f) \otimes \Omega \longrightarrow E \otimes \Omega
$$
is the natural
inclusion.
In this way, 
Lemma~\ref{prop-Phi-1-is-generically-surjective} provides 
the following alternative description of 
the moduli space $\Bpq$ of binary Hodge bundles:
{\em
$\Bpq$ parameterizes a family of equivalence classes of pairs 
$(E,f:E \rightarrow F)$ 
of a rank $p$ vector bundle $E$ of degree $d_Q$ and a length 
$p(2g-2)+d_Q-d_P$ quotient sheaf $f: E \rightarrow F$ on $X$. 
}
The irreducibility of $\Bpq$ is an easy consequence of this description.


\begin{pf} (of Proposition~\ref{prop:irreducible})
We construct an auxiliary irreducible variety $Q_2$ and 
a Zariski open subset $Q_2^{ss} \subset Q_2$.  The subset $Q_2^{ss}$
maps surjectively onto 
$\Bpq$.  The scheme
$Q_2$ is a relative Quot scheme over a subset $R$ of a Quot scheme $Q_1$
of vector bundles. We recall the construction of $R$ following
Seshadri \cite{Se1}.

The space $\Bpq$ is a subset of the moduli space of Higgs pairs. 
The family of Higgs bundles parameterized by $\C \M$ is bounded \cite{Si2}.
Hence, there exists an ample line bundle $L$ on $X$
with the following property:
{\em Every Hodge bundle in the family of Hodge bundles
parameterized by $\Bpq$ admits a
representative $(V_P\oplus V_Q,\Phi_1)$
with vanishing first cohomology $\H^1(X,V_Q\otimes L)=0$}.
Let $H(m):= \chi(V_Q\otimes L^{m})$ be the Hilbert polynomial
of rank $p$ vector bundles of degree $d_Q$. Set $a=H(1)$ and
let $Q_1:=Quot^H_{\oplus_{i=1}^a L^{-1}/X/\C}$ be the
Grothendieck scheme parameterizing
the quotient sheaves of $\oplus_{i=1}^a L^{-1}$ with Hilbert polynomial $H$ \cite{Gr1}.
The scheme $Q_1$ contains an irreducible and smooth
quasi-projective variety $R$ defined by
$$
R = \{W \in Q_1 : W \mbox{ is locally free and } \H^1(W) = 0 \}
$$
(\cite{Se1} Chapter III Proposition 23).
By our choice of $L$, every Higgs pair in
$\Bpq$ is represented by a pair $(V_Q\oplus V_P,\Phi)$,
such that $V_Q$ is realized in $R$ as a quotient of $\oplus_{i=1}^a L^{-1}$.

Let
$$
E \longrightarrow X \times R
$$
be the universal quotient bundle of $\oplus_{i=1}^a L^{-1}$.
Then there exists a relative Quot scheme
\[
Q_2 \ := \ Quot^{-d_P + d_Q + 2p(g-1)}_{E/X \times R/R}
\]
parameterizing quotient sheaves of $E$ supported as
length $-d_P + d_Q + 2p(g-1)$ subschemes of a fiber of
$X \times R \rightarrow R$ \cite{Gr1}.
By construction and Lemma~\ref{prop-Phi-1-is-generically-surjective},
$Q_2$ parameterizes a family of Higgs bundles that contains
representatives of all classes in $\Bpq$.

The morphism $Q_2\rightarrow R$ factors through
a surjective morphism
\[
h \ : \ Q_2 \ \ \rightarrow \ \
R\times X^{(-d_P + d_Q + 2p(g-1))},
\]
where $X^{(m)}$ stands for the $m$-symmetric product of $X$.
A quotient sheaf $F$ is sent to ${\displaystyle \sum_{x\in X} \ell_x\cdot x}$,
where $\ell_x$ is the length of the stalk of $F$ at the point $x$.
Each fiber of $h$ is isomorphic to the product of infinitesimal Quot schemes
$Q(\ell,\O_{(x)}^p)$ of length $\ell$ quotients of the stalk at $x$ of the
trivial rank $p$ vector bundle.
Consider the subscheme $Q_2^{free}$ of $Q_2$ parameterizing
pairs of quotient sheaves
$(\oplus_{i=1}^a L^{-1}\rightarrow W, W\rightarrow F)$, where $F$
is supported on a subscheme $D\subset X$ as a free $\O_D$-module of rank $1$.
The restriction of $h$ to $Q_2^{free}$
is a smooth morphism.
The scheme $Q_2^{free}$ is irreducible because $R$, 
the symmetric product of $X$,
and each fiber of $h$ are irreducible.
$Q_2^{free}$ is dense in $Q_2$,
because each fiber of $h$ in $Q_2^{free}$ is dense in the fiber
in $Q_2$ (Lemma \ref{lemma-irreducibility-of-a-fiber}).

Finally, the semi-stable condition is open.  Hence, the subscheme $Q_2^{ss}$ of
$Q_2$, parameterizing the semi-stable Higgs bundles, is open and, consequently,
irreducible.  Proposition~\ref{prop:irreducible} follows from the fact that
$Q_2^{ss}$ admits a surjective morphism onto $\Bpq$.
\end{pf}

\begin{lemma}
\label{lemma-irreducibility-of-a-fiber}
The Quot scheme $Q(\ell,\O_{(x)}^p)$, of length $\ell$ quotients of the stalk
at $x$ of the trivial rank $p$ vector bundle, is irreducible.
\end{lemma}
\begin{pf}
The Lemma must be well known.
We could not find a reference, so we include a short proof.
Let $m$ be the maximal ideal of $x$
and $A$ the ring $\O_{(x)}/m^\ell$.
Any quotient sheaf in $Q(\ell,\O_{(x)}^p)$ is also a quotient of
the free $A$-module of rank $p$. It is a direct sum of at most $p$
cyclic $A$-modules (Nakayama's Lemma and the classification of modules over
the discrete valuation ring $\O_{(x)}$).
Thus, the isomorphism class of a quotient sheaf
is determined by a partition of $\ell$ as a sum of $p$ non-negative integers.

Let $G:=GL(p,A)$ be the group of
automorphisms of the free $A$-module of rank $p$.
The group $G$ acts on $Q(\ell,\O_{(x)}^p)$.
Let $Q(\ell,\O_{(x)}^p)^{(\ell)}$ be the orbit of quotient sheaves that are
free $A$-modules of rank $1$. It suffices to prove that the orbit
$Q(\ell,\O_{(x)}^p)^{(\ell)}$ is dense in $Q(\ell,\O_{(x)}^p)$.

The proof is by induction on $p$. If $p=1$, then
$Q(\ell,\O_{(x)}^1)=Q(\ell,\O_{(x)}^1)^{(\ell)}$ is a point.
Let $(\ell_1,\ell_2, \dots, \ell_p)$, $\ell_i\geq 0$,
be a partition of $\ell$. We show that the orbit
$Q(\ell,\O_{(x)}^p)^{(\ell_1,\ell_2, \dots, \ell_p)}$
is in the closure of $Q(\ell,\O_{(x)}^p)^{(\ell)}$.
If some of the $\ell_i$ are zero, the proof reduced to the case of a smaller
$p$. Assume that none of the $\ell_i$ vanishes. Let $z$ be a local parameter
and consider the linear combination
$\eta(t):=\psi+t\varphi$, $t\in \C$, of the two $p\times p$ matrices
\begin{eqnarray*}
\psi & = &
\left(
\begin{array}{ccccccccccc}
0          &         &    & \cdots & 0 & z^{\ell_p}
\\
z^{\ell_1} &    0    &    & \cdots &  & 0
\\
0       & z^{\ell_2} & 0  & \cdots &  & 0
\\
\vdots  &         &    &        &  & \vdots
\\
0       &  \cdots &    &    0   & z^{\ell_{p-1}} & 0
\end{array}
\right) \ \ \ \ \ \mbox{and}
\\
\varphi & = &
\left(
\begin{array}{ccccccccccc}
1       & 0       & & \cdots &             & 0           & 0
\\
0       & 1       & &        &             & 0           & 
\\
        & 0       & &        &             & \vdots      & \vdots
\\
\vdots  & \vdots  & &        & 1           & 0           & 
\\
        & 0       & & \cdots & 0           & 1           & 0
\\
0       & 0       & & \cdots & 0           & 0           & z^\ell
\end{array}
\right).
\end{eqnarray*}
When $t=0$, $\eta(0)$ is equal to $\psi$ and its cokernel is a quotient
sheaf in the orbit $Q(\ell,\O_{(x)}^p)^{(\ell_1,\ell_2, \dots, \ell_p)}$.
If $t\neq 0$, then row reduction leads to the expression of the column
$\eta(t)_p$ as a linear combination of the first $p-1$ columns
(mod $z^\ell$):
\[
\eta(t)_p \ = \
\frac{z^{\ell_p}}{t}\eta(t)_1 -
\sum_{k=2}^{p-1}\left[
\left(\frac{-1}{t}\right)^k
z^{(\ell_p+\sum_{i=1}^{k-1}\ell_i)}
\right]
\eta(t)_k
+ \left[t+\left(\frac{-1}{t}\right)^{p-1}
\right] \cdot z^{\ell} \cdot e_p,
\]
where $\{e_1, \dots, e_p\}$ is the standard basis of $\O_{(x)}^p$.
Consequently,  if $t\neq 0$ and $t^p\neq (-1)^p$,
the cokernel of $\eta(t)$ is in
$Q(\ell,\O_{(x)}^p)^{(\ell)}$.
\end{pf}

\end{document}